\documentclass{article}
\begin{document}

Tsemo Aristide

College Boreal

1, Yonge Street, Toronto, ON

tsemo58@yahoo.ca

\bigskip
\bigskip

\centerline{\bf Some properties of schemes in groups theory and Top couples.}

\bigskip
\bigskip

\centerline{\bf Abstract.}

\bigskip

Let $G$ be a group,  and $H$ a $G$-group defined by an imbedding map $G\rightarrow H$; in [12] we have defined
  a topology on a subset of normal subgroups of $H$, the so-called prime ideals. In this work, we generalize this topology to other classes of groups.
 We hope that such topologies define new effective tools to tackle well-known problems in groups theory.
  We study properties of the topology defined in [12], for example we define the notions of functor of points, dimension and
  outline a Galois theory.

\bigskip
\bigskip

{\bf 1. Introduction.}

\bigskip
\bigskip

Classical algebraic geometry is the study of solutions of polynomial equations defined
on commutative rings. The nature of the objects involved here can be very different when one allows
the base ring to vary; for example varieties defined over the fields of complex numbers are complex manifolds
with singularities, thus,  inherit a separated topology from their complex structure, this is not the case of a variety defined over a finite field.
The first tools developed in classical algebraic geometry where related to topology or complex geometry
on the purpose to study complex schemes. In fact, Serre has shown in the G.A.G.A [8] deep relations between the
objects of these theories.

To study the variety defined on more general rings, some tools  issued from complex algebraic geometry have been adapted.
The most efficient work in that direction has been done by A. Grothendieck with the theory of schemes;
he has introduced topological concepts to objects which does not look a priori topologically by generalizing the Zariski topology. This approach has
allowed to give new interpretations of many conjectures and often led to their solutions.
We can for examples quote the Weil conjectures and more recently the Milnor conjecture.

\smallskip

Many generalizations of the theory of schemes have been achieved for non commutative algebras.
In group theory, many classical sets are defined by algebraic-like equations. Thus formally,
they are similar to the objects studied in classical algebraic geometry.
Recently, Baumslag, Myasnikov and Remeslennikov [2] have initiated the study of varieties of groups.
They consider a group $G$, and $G[X_1,...,X_n]$ the free product of $G$ and the free group
generated by a set of cardinality $n$. For any group $H$ endowed with an injective morphism
$G\rightarrow H$, the elements of $G[X_1,...,X_n]$ define polynomial functions on $H^n$ which took
their values in $H$. The solutions of these equations  define a Zariski-like topology.
These authors have adapted tools from the theory of algebraic varieties like
Zariski topology and irreducible subsets in groups theory.

In our paper [12], we have generalized the construction of Baumslag, Myasnikov and Remeslennikov by
defining a theory of schemes in this context,  The purpose of this work is to
continue to study this notion; we also generalize it by extending it to others subcategories
of the category of groups. We hope that these new concepts will also lead to the solutions of
some problems in groups theory by using the topological intuition.

\bigskip

\centerline{\bf Plan.}

1. Introduction

2. Some properties of groups

3. Topology in groups theory

4. The geometry of $G$-groups

5. The functor of points

6. Noetherian $G$-schemes

7. Dimension of affine schemes

8. Galois theory.

\bigskip

{\bf 2. Some properties of groups.}

\bigskip

Let $G$ be a group,  the comma category  over $G$ is the category
whose objects $(H,\phi_H)$ are  morphisms of groups $\phi_H: G\rightarrow H$.
 A morphism $f:(H,\phi_H)\rightarrow (H',\phi_{H'})$ of the comma category is  a morphism
of groups $f:H\rightarrow H'$ such that $f\circ \phi_H = \phi_{H'}$.

We denote by $C(G)$ the subcategory of the comma category such that each object
$(H,\phi_H)$ of $C(G)$ is an object of the comma category such that $\phi_H$ is injective.

\medskip

Let $(H,\phi_H)$ be an object of $C(G)$ and $x$ an element of $H$.
For every element $g\in G$, we denote $\phi_H(g)x\phi_H(g^{-1})$ by $x^g$. We denote by $G(x)$ the  subgroup
of $H$ generated by $\{x^g, g\in G\}$.
We say that the element $x\in H$ is invertible if and only if $G(x)\cap \phi_H(G)\neq 1$.

Let $H$ be a group, and $g$, $g'$ elements of $H$; the commutator $[g,g']$ of $g$ and $g'$ is the element
$gg'g^{-1}{g'}^{-1}$. If $U$ and $V$ are subgroups of $H$, $[U,V]$ is the normal  subgroup generated by
$[u,v], u\in U, v\in V$.

Let $H$ be a group, we denote by $D^0H$ the group $H$, suppose defined $D^nH$, for an integer $n$,
then $D^{n+1}H = [D^nH,D^nH]$.

We say that an element $x$ of $H$ distinct of the neutral element is a divisor of zero if and only if there exists an element
 $y$ of $H$ distinct of the neutral element
 such that the group of commutators $[G(x),G(y)]=1$. Remark that the element $y$
in this situation is also divisor of zero.

\medskip

{\bf Proposition 2.1.}

{\it The element $x$ of $H$ is a divisor of zero if and only if there exists a non trivial element $y$ of $H$
such that for every element $g\in G$,  $[x,y^g]=1$.}

\medskip

{\bf Proof.}

Let $x$ be an element of $H$ which is a divisor of zero, and $y$ an element of $H$ such that $[G(x),G(y)]=1$.
The set $U= \{[x,y^g],g\in G\}$ is contained in $[G(x),G(y)]$; we deduce that all  elements of $U$ are equal to the
neutral element.

Conversely, suppose that all elements of $U$ are equal to the neutral element. Let $u,v\in G$,
$$
[x^u,y^v] = uxu^{-1}vyv^{-1}ux^{-1}u^{-1}vy^{-1}v^{-1}.
$$

We deduce that:

$$
u^{-1}[x^u,y^v]u = xu^{-1}vyv^{-1}ux^{-1}u^{-1}vy^{-1}v^{-1}u =[x,y^{u^{-1}v}]
$$

Since $[x,y^{u^{-1}v}] = 1$,we deduce that $[x^u,y^v]=1$.

Recall that if $a,b,c$ are elements of a group, if $a$ commutes with $b$ and $c$, then $a$ commutes with $b^{-1}$ and
$bc$, this implies that  $G(x)$ commutes with $G(y)$.

\medskip

{\bf Remark.}

\medskip

Let $N_G(x)$ be the normal subgroup generated by $G(x)$, we can find an example of a $G$-group  $H$ such that
 $[G(x),G(y)] =1$ and $[N_G(x),N_G(y)]\neq 1$:

Consider $H=SL(4,C)$ the group of linear automorphisms of the $4$-dimensional complex vector space $C^4$ whose determinant is $1$.
Let $e_1,e_2,e_3,e_4$ be a basis of $C^4$. Let $G = SL(2,C)$ be the subgroup of $H$ such
that for every $f\in G$, $f(e_3) =e_3, f(e_4) =e_4$ and $f(e_1)$ and $f(e_2)$ are elements of the vector subspace $Vect(e_1,e_2)$
generated by $e_1$ and $e_2$. Let $g$ be an element of $H$ such that $g(e_1) =e_1, g(e_2) = e_2$
and $g(Vect(e_3,e_4))\subset Vect(e_3,e_4)$, $[g,G] =1$, but $[N_G(g),G]$ is not the trivial subgroup: if $f:C^4\rightarrow C^4$
is defined by $f(e_1) = e_3, f(e_2) = e_4, f(e_3) =e_1$ and $f(e_4) =e_2$, $fgf^{-1}\in G$, and we can choose $g$ such
that $fgf^{-1}$ is not in the center of $G$.

\medskip

{\bf Proposition 2.2.}

{\it Suppose that $G$ is a   solvable $G$-group different of $\{1\}$, then any $G$-group $H$ has a divisor of zero. If $H$ is a finite dimensional simply connected Lie group
which does not have zero divisors, then $H$ is simple.}

\medskip

{\bf Proof.}

Let $G$ be a   solvable group different of $\{1\}$;  there exists $n$ such that $D^nG$ is not trivial and $D^{n+1}G$ is trivial.
For every element $x\in D^nG$ different of
its neutral element we have $[G(x), G(x)] =1$, thus $x$ is a divisor of zero; if
$H$ is a $G$-group, $x$ is also a zero divisor of $H$.

Suppose that $H$ is a finite dimensional simply connected Lie group which does not have divisors of zero and its radical is not trivial.
  let $R$ be the radical of $H$, recall that $R$ is a maximal solvable
normal subgroup of $H$. Let $n$ be an integer such that $D^nR$ is not trivial, and $D^{n+1}R$ is trivial.
Since $R$ is a normal subgroup of $H$, for every element $x\in D^nR$ and every element $g\in G$, $x^g\in D^nR$,
this implies that $[G(x),G(x)] =1$ for every element $x\in D^nR$. In particular, if $x$ is not the neutral element
it is a divisor of zero. This is in contradiction with the fact that we have supposed that $G$ does not
have divisors of zero.
Thus the radical of $H$ is trivial and $H$ is semi-simple. If $H$ is not simple, then it is the product of
two non trivial groups $H_1$ and $H_2$. Let $h_1$ and $h_2$ be respectively
elements of $H_1$ and $H_2$ distinct of the neutral element, we have
$[G((h_1,1));G((1,h_2))] =1$. Thus $H$ has divisors of zero.

\medskip

{\bf Remark.}

\medskip

The previous proposition shows that a non commutative nilpotent group $G$ has divisors of zero.
In fact, every element in the center of $G$ is a divisor of zero.

\medskip

{\bf Proposition 2.3.}

{\it Let $k$ be a field, and $G$  a subgroup of $PSL(n,k)$ suppose that the action of $G$ on the $n$-dimensional
$k$-projective space $Pk^n$ is reducible, then $G$ has divisors of zero.
If $G$ is a connected semi-simple Lie group defined over the complex numbers imbedded in $PSL(n,C)$ or if $G$ is a finite group
imbedded in $PSL(n,k)$ where $k$ is a field whose characteristic does not divide the order of $G$,
then $G$ does not preserves  the image of a proper linear vector space in $P^nC$ (resp. $P^nk)$
if $G$ does not have zero divisors.}

\medskip

{\bf Proof.}

Suppose that the action of $G$ is reducible. Let $G'$  be the inverse image of the group $G$ by the projection map $p:Sl(n,k)\rightarrow PSL(n,k)$.
 We can write $k^n = U\oplus V$, where $U$ and $V$ are non trivial vector subspaces of $k^n$ and are stable by $G'$.
Let $u'$ be a non trivial invertible linear map of $U$ and $v'$ be a non trivial invertible linear map of $V$. Consider the linear morphism
$u$ of $k^n$ defined by $u(x)=u'(x)$ if $x\in U$, $u(x) = x$ if $x\in V$. Let $v$ be the linear morphism of $k^n$
defined by $v(x) = x$ if $x\in U$, $v(x) = v'(x)$ if $x\in V$.  We have $[G'(u),G'(v)]=1$. It results that
the image of $u$ and $v$ in $PSL(n,k)$ are divisors of zero.

Suppose that $G$ is a connected semi-simple  subgroup of $PSL(n,C)$, (resp. $G$ is a finite subgroup of $PSl(n,k)$
such that the order of $G$ and the characteristic of $k$ are prime each other),
$p^{-1}(G)$ is semi-simple (resp. finite subgroup) and $C^n$ (resp $k^n$) under the action of $G$ can be decomposed
in the direct summand of irreducible vector spaces. Thus $PSL(n,C)$ (resp. $PSL(n,k)$ has zero divisors
if  $p^{-1}(G)$ preserves a proper linear subspace.

\medskip

{\bf Definitions 2.1.}

A $G$-domain $(H,\phi_H)$ of $C(G)$ is an object of $C(G)$ such that $H$ does not have divisors of zero.

A normal subgroup  $I$ of $H$ is an ideal of $H$ if and only if $I\cap G = \{1\}$.

A normal subgroup $P$ of $H$ is a prime ideal if and only if:

- $P$ is an ideal distinct of $H$ and

- $G/P$ is a $G$-domain, this fact is equivalent to saying that for every $x,y\in H$, $[G(x),G(y)]\in P$
if and only if $x$ is in $P$ or $y$ is in $P$.

- The $G$-group $H$ is $G$-simple if $\{1\}$  is the  the only ideal of $H$.

\medskip

{\bf Proposition 2.4.}

{\it Suppose that $G$ is a connected simple Lie group over a field of characteristic zero, and $H$ a  connected Lie group
in $C(G)$,  if $H$ is $G$-simple  then $H$ is simple Lie group.}

\medskip

{\bf Proof.}

We can write $H$ has a semi-product of a semi-simple group $U$ and a solvable group $L$ by using the Levi
decomposition. Without restricting the generality, we can suppose that $G$ is contained in $U$ since $G$ is simple.
Since $G\cap L = \{1\}$, $L$ is an ideal. We deduce that $L$ is trivial since the only prime of $H$ is $1$.
If the subgroup $U$ is not simple,  there will exist a simple component of $U$ which intersection  with $G$ is trivial since
the simple Lie group $G$ is contained in a simple component of $H$. This fact is
impossible since all ideals of $H$ are trivial.

\medskip

{\bf Remark.}

\medskip

The previous result is not true if we don't assume that the groups $G$ and $H$ are connected as shows
the following example: Let $S_n$ be the group of permutations of the set of $n$ elements, $n>4$. The
subgroup $A_n$ of permutations of even signature is a simple group. The group $S_n$ is endowed with an $A_n$-structure.
It is well-known that the only normal subgroups of $S_n$ are the trivial subgroup, $A_n$ and $S_n$. Thus
$S_n$ is an $A_n$-simple group, and $S_n$ is not simple.

\medskip

{\bf Proposition 2.5.}

{\it Let $k$ be a field and $G$ a Zariski dense subgroup of $PSl(n,k)$,  the $G$-group $PSl(n,k)$ does
not have divisors of zero.}

\medskip

{\bf Proof.}

Let $x$ be a divisor of zero and $y$ a non trivial element of $PSl(n,k)$ such that for every element $g\in G$,
$[x,y^g]=1$. The equation $[x,y^u]=1$ is an algebraic equation in $PSl(n,k)$. This implies that the Zariski adherence
of $G$ satisfies this equation. Since $G$ is Zariski dense in $PSL(n,k)$, this implies that for every $g\in PSl(n,k)$,
$[x,y^g] = 1$. Let $x'$ (resp. $y'$) an element of $Sl(n,k)$ above $x$ (resp. above $y$). Thus $[x',{y'}^g] = cI_n$
is in the center of $SL(n,k)$ for every $g\in SL(n,k)$. Recall that the elements of the center of $SL(n,k)$ are $c_1I_1,...,c_nI_n$ where $c_1,...,c_n$ are
the $n$-roots of unity.
  The subset $G_i$ of $SL(n,k)$ such that $[x',{y'}^g] = c_i^{-1}I_n$ is closed, and  $SL(n,k)$
is a finite union of the $G_i$. Since $SL(n,k)$ is irreducible, there exists a $i$ such that $G_i = SL(n,k)$.  It
results that $[{x'}^n,{y'}^g] =I_n$ for every $g\in SL(n,k)$. Let $E_h$ be an eigenvector space of $y'$
associated to the eigenvalue $h$. Let $p$ be the dimension of $E_h$, for every vector subspace $E$ of dimension $p$,
 there exists an element of $SL(n,k)$ such that $g(E_h) = E$, this implies that $E$ is the eigenvector space of ${y'}^g$
 and ${x'}^n$ preserves $E$. We deduce that ${x'}^n =cI$, $c\in k$. The same argument shows also that ${y'}^n = dI_n$, $d\in k$.
 Thus $x'$ and $y'$ are diagonalizable maps.  If $g$ is an element of $Sl(n,k)$, $g(E_h)$ is an eigenvector space of
${y'}^g$ associated to $h$. This implies that ${x'}(g(E))$ is the eigenvector space  of ${y'}^g$ associated to $c_i^{-1}h$ if
$[x',{y'}^g] = c_iI_n$. Thus $x'$ permutes the eigenvector spaces of $y'$. Since $x'$ and $y'$ are not homothetic maps since $x$
and $y$ are not trivial,  $y'$ has two
 distinct eigenvector spaces $E_1$ and $E_2$ associated respectively to the eigenvalue $c_1$ and $c_2$;
  without restricting the generality, since $x'$ permutes the eigenvector space of ${y'}^g$ for every $g\in SL(n,k)$, we can suppose that
   $x'(g(E_1)) = g(E_2)$. This last equality is not
 possible for every element of $SL(n,k)$.

\medskip

{\bf Proposition 2.6.}

{\it Let $H$ be an element of $C(G)$, suppose that $G$ is a normal subgroup of $H$, then for every ideal
ideal $I$, $[I,G]=1$.}

\medskip

{\bf Proof.}

Let $g$ be an element of $G$, and $i$ and element of $I$, we have $[u,i] =uiu^{-1}i^{-1}$ is an element of
$G\cap I = \{1\}$.

\medskip

{\bf Definitions 2.2.}

Let $(H,\phi_H)$ be an object of $C(G)$, an element $x$ of $H$ is nilpotent of length $n$ if and only
if the group $G(n)(x)= [G(x),[G(x)[,...,G(x)] = 1$ where in the formula $G(x)$ appears $n$ times.

Let $H$ be an element of $C(G)$, the radical $rad_G(H)$ of $H$ is the intersection of all prime ideals of $H$.
We denote by $Com_H(G)$ the subgroup of $H$ such that for every element $h\in Com_H(G)$ and $g\in G$, we have
$[g,h] = 1$. The group $Com_H(G)$ is included in every element $P\in Spec_G(H)$ and so is its normalizer $N(Com_H(G))$.
This implies that $N(Com_H(G))$ is a subset of $Rad_G(H)$.

\medskip

{\bf Proposition 2.7.}

{\it Let $P$ be a prime and $x$ an element of $H$, suppose that $G(n)(x)$ is contained in $P$, then $x$ is
an element of $P$.  In particular The normal subgroup  $Nil(H)$ generated by nilpotent elements of $H$ is contained in $P$.}

\medskip

{\bf Proof.}

We show the result recursively.
 Let $P$ be a prime, suppose that every element $u$ such that $G(n-1)(u)\in P$ is in $P$.
 Let $x$ be an element of $H$ such that,   $G(n)(x)$ is contained in $P$.
Suppose that $x$ is not in  $P$. Let $u$ be an element of $G(n-1)(x)$,
we have $[G(x),u]\subset G(n)(x) \subset P$. Since $P$ is prime, $u$ is an element of $P$, we deduce that
$G(n-1)(x)$ is contained in $P$, the  iteration hypothesis implies that $x$ is an element of $P$.

\medskip

{\bf Proposition 2.8.}

{\it Suppose that $G$ is a simple group, and $H$ is an element of $C(G)$, then every normal subgroup of $H$
is an ideal or contain $G$.}

\medskip

{\bf Proof.}

Let $I$ be a normal subgroup of $H$, $I\cap G$ is a normal subgroup of $G$; since $G$ is simple we deduce that
$I\cap G = G$ or $I\cap G =\{1\}$.

\bigskip

{\bf 3. Topologies in group theory.}

\bigskip

In this part, we are going to generalize the topology introduced in [12]. Let $C$ be a subcategory  of the category of groups,
we denote by $D$ a subclass of the class of objects of $C$ which satisfies the following properties:

T1 If $G$ is an object of $C$, $G'$ is an object of $D$ such that there exists an injective $C$-morphism $i:G\rightarrow G'$,
then $G\in D$.

T2 Let $G$ be an element of $D$, and $I$, $J$ normal subgroups of $G$,  $[I,J] = 1$ implies that $I = 1$ or $J = 1$.

\medskip

{\bf Definitions 3.1.}

Let $H$ be an element of $C$, an ideal of $H$ is a normal subgroup $I$ of $H$ such that $H/I\in C$ and the
projection $H\rightarrow H/I$ is a morphism of $C$.

An ideal $P$ of $H$ is a prime if and only if $H/P$ is an object of $D$ and $P$ is distinct of $H$.

T3 We suppose that the inverse image of an ideal by a morphism of $C$ is an ideal.

A couple $(C,D)$ which satisfies the properties T1, T2 and T3 will be called a Top couple.

\medskip

Let $H$ be an element of $C$ and $I$ a normal subgroup of $H$, we denote by $V(I)$ the set of prime ideals of $H$
which contain $I$.

\medskip

{\bf Proposition 3.1.}

{\it Let $H$ be an element of $C$, and $I$, $J$ normal subgroups of $H$, we have:

$$
V([I,J]) = V(I)\bigcup V(J)
$$

Let $(I_a)_{a\in A}$ a family of normal subgroups of $H$, and $I_A$ the normal subgroup generated by $(I_a)_{a\in A}$ we have:

$$
V(I_A) =\cap_{a\in A}V(I_a)
$$
}

\medskip

{\bf Proof.}

Let $P$ be an element of $V(I)\bigcup V(J)$, $P$ contains $I$ or $J$. This implies that $P$ contains $[I,J]$.

Let $P$ be an element of $V(I,J])$ and $p:H\rightarrow H/P$ the natural projection map, we have $[p(I),p(J)] = 1$.
Since $H/P$ is an object of $D$, we deduce that $p(I) =1$ or $p(J) = 1$. This is equivalent to saying that $I\subset P$ of $J\subset P$.

Let $P$ be an element of $V(I_A)$, for every $a\in A$, we have $I_a\subset I_A\subset P$ this implies $P\in \cap_{a\in A}V(I_a)$.

Let $P$ be an element of $\cap_{a\in A}V(I_a)$, for every $a\in A, I_a\subset P$ this implies that
$I_A\subset P$ and $P\in V(I_A)$, we deduce that $V(I_A) = \cap_{a\in A}V(I_a)$.

\medskip

{\bf Remarks.}

Let $(C,D)$ be a Top couple and $H$ be an object of $C$, and $Spec(H)$ the set of prime ideals of $H$. The set $V(I)$ are the closed subsets
of a topology defined on $Spec(H)$; $V(H)$ is the empty subset, and $V(\{1\}) = Spec(H)$.

\medskip

{\bf Proposition 3.2.}

{\it Let $(C,D)$ be a Top couple and $H$, $H'$ two objects of $C$. A $C$-morphism $f:H\rightarrow H'$ induces
a continuous morphism $c(f): Spec(H')\rightarrow Spec(H)$ by setting $c(f)(P) = f^{-1}(P)$.}

\medskip

{\bf Proof.}

Let $P$ be an ideal of $H'$, the property T3 implies that $f^{-1}(P)$ is an ideal subgroup of $H$
since the morphism $H/f^{-1}(P)\rightarrow H'/P$ is injective,
the property T1 implies that $H/f^{-1}(P)\in D$. This implies that $f^{-1}(P)$ is a prime ideal.

The morphism $c(f)$ is continuous since $c(f)^{-1}(V(I)) = V(f(I))$.

\medskip

{\bf Examples of Top couples.}

\medskip

We are going to present examples of Top couple, we start by the following example which is already presented in
Tsemo [12]:

\medskip

{\bf Proposition 3.3.}

{\it Let $G$ be a group, and $D$ the class of $G$-groups which  do not have divisors of zero, the couple $(C(G),D)$ is a Top couple.}

\medskip

{\bf proof.}

We have to show that the couple $(C(G),D)$ satisfies the properties $T1$, $T2$ and $T3$.

If $H$ is a $G$-group which does not have zero divisors, any $G$-subgroup of $H$ does not have zero divisors of zero, thus
the property $T1$ is satisfied.

We show now that the couple $(C(G),D)$ satisfies the property T2.
Let $H$ be a group which does not have divisors of zero, and $I,J$ two normal subgroups of $H$ such that
$[I,J] = 1$. Suppose that $I\neq 1$ and $J\neq 1$. This implies the existence of elements $x\in I$ and $y\in J$
different of the neutral elements. We have $[G(x),G(y)]\subset [I,J] = 1$. This is equivalent to saying that $x$ and $y$
are divisors of zero. This fact is a contradiction to the hypothesis; thus $I =1 $ or $J=1$.

It remains to show that the couple $(C,D)$ satisfies the property $T3$. Recall that an ideal $I$ of the element $H$
of $C(G)$ is a normal subgroup $I$ of $H$ such that $I\cap G =1$, if $f:H\rightarrow H'$ is a $G$-morphism, and
$I$ an ideal of $H'$, $f^{-1}(I)$ is an ideal of $H$ since the restriction of $f$ on $G$ is injective.

\medskip

We present now another example of a Top couple. One of the most interesting problem in groups theory is to determine
if a group has a finite subgroup which is isomorphic to the free group generated by two elements. Tits Jacques has shown
that a subgroup of a linear group is virtually solvable or he contains a free subgroup isomorphic to $F_2$.
The Top couple that we are presenting now is related to this subject.

Here $C$ is $C(F_2)$ the category whose objects are injective morphisms $F_2\rightarrow H$. An object $F_2\rightarrow H$ of
$C$ is an object of $D$ if and only if $H$ is a free group. Remark that the rank of an object of $D$ is greater or equal to $2$.
Let $(H,\phi_H)$ be an object
of $C(F_2)$, an ideal $I$ of $H$ is a normal subgroup $I$ such that $I\cap \phi_H(F_2) = \{1\}$. From this definition, we deduce that the
ideal $P$ of $H$ is a prime if and only if $H/P$ is a free group.

\medskip

{\bf Proposition 3.4.}

{\it The couple $(C(F_2),D)$ is a Top couple.}

\medskip

{\bf Proof.}

Let show that $(C(F_2),D)$ satisfies $T1$. Let $i:H\rightarrow H'$ be an injective $F_2$-morphism,
the fact that $H'$ is an object of $D$ is equivalent to saying that $H'$ is free, this implies that
$H$ is free since a subgroup of a free group is free.

Property T2. Let $H$ be an object of $D$ and $I,J$ ideal of $H$ such that $[I,J] = 1$. Since $H$ is a free group, we
deduce that $I =1$ or $J=1$.

Property T3 Let $f:H\rightarrow H'$ be an $F_2$-morphism, and $P$ an ideal of $H'$, $f^{-1}(P)$ is an ideal of $H$
since the restriction of $f$ on $F_2$ is injective, the canonical map $H/f^{-1}(P)\rightarrow H'/P$ is injective
and a subgroup of a free subgroup is free.

\bigskip

{\bf 4. The geometry of $G$-groups.}

\medskip

In the sequel we are going to consider only the topology mentioned at proposition 3.3.

\medskip

{\bf Proposition 4.1.}

{\it Let $H$ be an element of $C(H)$, for every normal subgroups $I$ and $J$ of $H$, we have $V([I,J])= V(I\cap J)$.}

{\bf Proof.}

We know that $V([I,J]) = V(I)\bigcup V(J)$. Since $[I,J]\subset I \cap J$, we have $V(I\cap J)\subset V([I,J])$.
Let $P$ be an element of $V([I,J])$, suppose that $P$ is not an element of $V(I\cap J)$, this implies that there exists
an element $x$ of $I\cap J$ which is not an element of $P$. For every element $g\in G$, $[x,x^g]\in [I,J]$ since $x\in I\cap J$.
Since $P$ is a $G$-domain, we deduce that $x$ is an element of $P$. This is a contradiction, thus $P$ contains $I\cap J$.

\medskip

{\bf Proposition 4.2.}

{\it Suppose that $G$ is a solvable group distinct of $\{1\}$, for every object $H$ of $C(G)$, $Spec_G(H)$ is empty.}

\medskip

{\bf Proof.}

Let  $P$ be a prime ideal of an object $H$ of $V(I)$. Since $G$ is a solvable group distinct of $1$, the proof of the proposition 2.2
shows that there exists an element $x\in G$ which is a divisor of zero. This implies that $x\in P$. This is impossible
since $P\cap G = \{1\}$.

\medskip

{\bf Remark.}

\medskip

Let $G$ be the trivial subgroup $\{1\}$, any group $H$ is an element of $C(G)$. Let $P$ be a prime ideal of $H$,
for any element $x\in H$, $x\in P$ since $[x,x] =1$. This implies that $Spec_G(H)$ is empty since a prime of $H$ is distinct of $H$.

Let $H$ be an element of $C(G)$, the proposition 3.3 shows that the subsets $V(I)$ where $I$ is a normal subgroup of $H$ are the closed subsets
of a topology on the set $Spec_G(H)$ of prime ideals of $H$.  The empty closed subset is $V(H)$, and the $Spec_G(H) = V(\{1\})$.

We will say that $Spec_G(H)$ is an affine scheme.

\bigskip

{\bf The structural sheaf.}

\medskip

Let $(H,\phi_H)$ be an element of $C(G)$, and $U$ an open subset of $Spec_G(H)$, we denote by $Rad(U)$ the
intersection of the elements of $U$. We define on $Spec_G(U)$ the structural presheaf:

$$
U\longrightarrow P_H(U) = H/Rad(U)
$$

Let $U$ and $V$ be two open subsets of $Spec_G(H)$ such that $U\subset V$; $Rad(V)\subset Rad(U)$, this implies
the existence of a map $r_{U,V}:P_H(V)=H/Rad(V)\rightarrow P_H(U) = H/Rad(U)$ which is the restriction map of the
presheaf $P_H$.

The structural sheaf of $O_H$ of $Spec_G(H)$ is the sheaf associated to the presheaf $P_H$.

Recall that if $x$ is an element of $Spec_G(H)$ and $U(x)$ the set of open subsets of $Spec_G(H)$ which contains $x$,
 the stalk ${P_H}_x$ of $x$ is the inductive limit ${P_H}_x = lim_{V\in U(x)}P_H(V)$. The etale space $E_{P_H}$ of $P_H$ is the topological space
 which is the disjoint union $\{{P_H}_x, x\in Spec_G(H)\}$ endowed with the topology generated by $\{ (s_x)_{x\in U}, s\in P_H(U)\}$
 where $U$ is any open subset of $Spec_G(H)$. The sheaf $O_H$ is a the sheaf of continuous sections of the canonical
 morphism $p:E_{P_H}\rightarrow Spec_G(H)$ which sends an element of ${P_H}_x$ to $x$.

\medskip

{\bf Definition 4.1.}

Let $(X,O_X)$ and $(Y,O_Y)$ be two topological spaces $X$ and $Y$ respectively endowed with the
sheaves $O_X$ and $O_Y$. A morphism $(f,f'):(X,O_X)\rightarrow (Y,O_Y)$ is defined by
 a continuous map $f:X\rightarrow Y$ and a morphism of sheaves $f':O_Y\rightarrow f^{-1}O_X$ which commutes
 with the restriction maps.

 \medskip

{\bf Proposition 4.3.}

{\it Let $(H,\phi_H)$ and $(L,\phi_L)$ elements of $C(G)$ with trivial radicals, there exists a natural bijection between morphisms:
 $(f,f'):(Spec_G(H),O_H)\rightarrow (Spec_G(L),O_L)$ and $Hom_G(L,H)$.}

 \medskip

 {\bf Proof.}

 Let $(f,f'):(Spec_G(H),O_H)\rightarrow Spec_G(L),O_L)$ be a morphism. We associate to $(f,f')$ the $G$-morphism
 $f'(Spec_G(L)):O_L(Spec_G(L))=L \rightarrow O_H(Spec_G(H))=H$.

 Let $u:L\rightarrow H$ be a $G$-morphism of $C(G)$, then $u$ induces a continuous map
 $(f_u,f'_u):(Spec_G(H),\phi_H)\rightarrow (Spec_G(L),\phi_L)$
defined as follows:

Let $P$ be an element of $Spec_G(H)$, then $f_u(P)=u^{-1}(P)$. The proposition 3.3 shows that
 $f_u$ is well defined.  Let $V$ be an open subset of $Spec_G(L)$, the composition by $u$ induces a morphism
 $f'_u(V):O_{L}(V)\rightarrow f^{-1}(O_H)(V)$.

The correspondences defined above between $Hom_G(L,H)$ and the set of morphisms between  $(Spec(H),O_H)$ and $(Spec(L),O_L)$ are inverse each others.

\medskip

{\bf Definition 4.2.}

A topological manifold $X$ is an $G$-scheme if and only if there exists a covering $(U_i)_{i\in I}$ of $X$
and a sheaf $O_X$ on $X$, such that for every $i\in I$, there exists an affine $G$-scheme $(Spec_{G}(H_i),O_{H_i})$
and an homeomorphism of sheafed spaces $f_i:(Spec_{G}(H_i),O_{H_i})\rightarrow (U_i,O_{U_i})$, where $O_{U_i}$
is the restriction of $O_X$ to $U_i$.

\medskip

{\bf Proposition 4.4.}

{\it Let $(H,\phi_H)$ be an element of $C(H)$, and $Nil(H)$ the normal subgroup of $H$ generated by its nilpotent elements.
The projection morphism $p':H\rightarrow H/Nil(H)$ induces
an homeomorphism between $p:Spec_G(H/Nil(H))\rightarrow Spec_G(H)$.}

\medskip

{\bf Proof.}

Let $P$ be a prime ideal of $H/Nil(H)$, $p'({p'}^{-1})(P)=P$. This implies that $p$ is injective.
Let $Q$ be a prime ideal of $H$, ${p'}^{-1}(p'(Q))=Q$ since $Q$ contains $Nil(H)$. This implies that $p$ is surjective thus bijective.
Let $I$ be an ideal of $H/Nil(H)$, $p(V(I))= V({p'}^{-1}(I))$, this implies that $p$ is open. We conclude  that $p$ is
an homeomorphism.

\bigskip

{\bf Definition 4.3.}

Let $G$ and $L$ be two elements of $C(G)$ with trivial radical. A morphism $(f,f'):(Spec_G(H),O_H)\rightarrow (Spec_G(L),O_L)$ is of finite type if and only if  the $L$-group
$H$ is the quotient of a free group $L[x_1,...,x_n]$ by an ideal generated by a finite subset
for the $L$-structure of $H$ induced by $(f,f')$.

 The smallest $n$ such that $H$ is the
quotient of $L[x_1,...,x_n]$ by an ideal of $L[x_1,...,x_n]$ generated by a finite subset is called the
$L$-rank of $G$.

Suppose that $G$ is a domain and the radical of $H$ is trivial, we say that $(Spec_G(H),)O_H)$ is of finite type if and only if the canonical  morphism
$(Spec_G(H),O_H)\rightarrow (Spec_G(G),O_G)$ is a morphism of finite type.

\medskip

We endow now the category of affine schemes with a product:

\medskip

{\bf Proposition 4.5.}

{\it Let $H$ and $K$ be two elements of $C(G)$, the affine scheme $Spec_G(H*_GK)$ is the categorical product of
$Spec_G(H)$ and $Spec_G(K)$.}

\medskip

{\bf Proof.}

A couple of morphisms $f:Spec_G(L)\rightarrow Spec_G(H)$ and $g:Spec_G(L)\rightarrow Spec_G(K)$ is defined by a couple of morphisms
of $G$-groups $f':H\rightarrow L$ and $g':K\rightarrow L$ which induces a unique morphism $(f',g'):H*_GK\rightarrow L$
such that $f' = (f',g')\circ i_H$ and $g' = (f',g')\circ i_K$ where $i_H:H\rightarrow H*_GK$ is the map defined by
$i_H(h) = h*_G1$. This is equivalent to saying that $Spec_G(H*_GK)$ is the product of $Spec_G(H)$ and $Spec_G(K)$.

\medskip

{\bf Proposition 4.6.}

{\it The category of affine $G$-schemes has sums.}

\medskip

{\bf Proof.}

Let $H$ and $K$ be two elements of $C(G)$, we can endow the product of these groups $H\times K$ with the diagonal
action of $G$; Let $f:Spec_G(H)\rightarrow Spec_G(L)$ and $g:Spec_G(K)\rightarrow Spec_G(L)$ be morphisms of $G$-schemes.
These morphisms defined respectively by $G$-morphisms: $f':L\rightarrow H$ and $g':L\rightarrow K$ which induce
a $G$-morphism $(f',g'):L\rightarrow H\times K$, which defines the sum $(f,g):Spec_G(H\times K)\rightarrow Spec_G(L)$.

\medskip

{\bf Definitions 4.4.}

Let $H$ be an element of $C(G)$, we define a functor $c_H:C(G)\rightarrow C(H)$ which associates to $L$ the amalgamated sum
$H*_GL$. We say that the $H$-scheme $Spec_H(H*_GL)$ is obtained from $L$ by a change of the basis. In particular if $H$ is
an algebraic closure $G_{al}$ of the group $G$, $Spec_{G_{al}}(G_{al}*_GL)$ will be called the geometric scheme associated
to $H$.

\medskip

{\bf Definitions 4.5.}

Recall that an irreducible topological set is a set which cannot be the union of two proper distinct closed subsets.

Let $(H,\phi_H)$ be an element of $C(G)$. The normal subgroup $I$ of $H$ is a radical, if $I$ is the intersection of
all the primes ideals which contains $I$.

\medskip

{\bf Proposition 4.7,}

{\it Let $(H,\phi_H)$ be an element of $C(H)$, if $I$ a prime ideal of $H$, then $V(I)$ is irreducible.
Conversely, suppose that $I$ is a radical ideal, then if $V(I)$ is irreducible $I$ is a prime ideal.}

\medskip

{\bf Proof.}

Suppose that $I$ is a prime ideal and $V(I)$ is the union of the proper closed subsets $V(J)$ and $V(K)$. Since
$I$ is a prime ideal, $I$ is an element of $V(I)$; thus $I$ is an element of $V(J)$ or $I$ is an element
of $V(K)$. If $I$ belongs to $V(J)$; this implies that $I$ contains $J$, thus $V(I)\subset V(J)$. This
is a contradiction with the fact that $V(J)$ is a proper closed subset of $V(I)$. The same argument shows
that $I$ cannot belong to $V(K)$.

Suppose now that $V(I)$ is an irreducible closed subset and $I$ is a radical ideal. Let $x$ and $y$ be elements of $H$ such that $[G(x),G(y)]\subset I$.
For every element $P\in V(I)$, $[G(x),G(y)]\subset P$, this implies that $x\in P$ or $y\in P$. Let
$V_x$ the subset of elements of $V(I)$ which contains $x$ and $V_y$ the subset of elements of $V(I)$
which contains $y$, $V(I)$ is the union of $V(x)$ and $V(y)$. Let $N(x)$ be the minimal subgroup of $H$ which contains $x$,
if $P$ is an ideal which contains $x$, $N(x)\subset P$. This implies that $V(x) = V(I)\cap V(N(x))$
and $V(y) = V(I)\cap V(N(y))$. Since $V(I)$ is irreducible, we deduce that $V(I)\cap V(N(x)) = V(I)$
or $V(I)\cap V(N(y)) = V(I)$. Suppose that $V(I)\cap N(x) = V(I)$, this implies that $V(I)\subset V(N(x))$
and for every prime ideal $P\in V(I)$, $x\in P$, since $I$ is a radical ideal, $x\in I$. By the same argument, we deduce that
if $V(I)\cap V(N(y)) = V(I)$, $y\in I$.

\medskip

{\bf Corollary 4.1.}

{\it Suppose that $G$ does not have divisors of zero, $Spec_G(H)$ is irreducible if and only if the radical of $H$ is a prime ideal.}

\medskip

{\bf Proof.}

Recall that the radical $Rad_G(H)$ is the intersection of all the prime ideals of $Spec_G(H)$. Thus by definition,
$Rad_G(H)$ is a radical ideal, we can apply proposition 4.7.

\medskip

{\bf Remark.}

If $Rad_G(H)$ is a prime, it is a generic point of $Spec_G(H)$, that is, its adherence is $Spec_G(H)$.

\bigskip

{\bf 5. Functor of points.}

\medskip

Let $H$ be an element of $C(G)$. Consider the $n$-uple $[h]=(h_1,...,h_n)$ of $H^n$; it induces a $G$-morphism:
$$
u_{[h]}:G[X_1,...,X_n]\rightarrow H
$$

defined by $u_{[h]}(X_i) =h_i$, for every element $f\in G[X_1,...,X_n]$, we set $f_H([h]) = u_{[h]}(f)$, we
have thus defined a morphism $f_H:H^n\rightarrow H$.

Let $S$ be a subset of $G[X_1,...,X_n]$ in their paper [2], Baumslag and al defined the subset:
$$
V_H(S) = \{[h] = (h_1,...,h_n)\in H^n: \forall f\in S, f_H([h]) = 1.\}
$$

We are going to present this notion with the concept of functor of points.

\medskip

{\bf Definition 5.1.}

We denote by  $Spec_G^0$ opposite to the category $Spec_G$ of affine $G$-schemes.
For every object $H$ of $C(G)$, the functor of points $h_{Spec_G(H)}: Spec_G^0\rightarrow Sets$ is the functor defined by:
$h_{Spec_G(H)}(Spec_G(K)) =mor_G(Spec_G(K),Spec_G(H)) = Mor_G(H,K)$.

\medskip

In particular, if $H$ is a scheme of finite type, $H$ is the quotient of $G[X_1,...,X_n]$ by an ideal $I$,
an element $f\in Mor_G(H,K)$ is defined by a $n$-uple $(k_1,...,k_n)\in K^n$ such that for every $f\in I$, $f(k_1,...,k_n) = 1$.
Thus the elements of $Mor_G(H,K)$ in this situation correspond to the elements of the closed subset $V_G(K)$ defined in [2], where
$S$ is replaced by $I$.

\bigskip

{\bf 6.1 Noetherian $G$-schemes.}

\medskip

The elements of $G[X_1,...,X_n]$ are often called equations in $n$-variables over $G$.
Every set of equations $E$ generates a normal subgroup $I_E$ of $G[X_1,...,X_n]$.

\medskip

{\bf Definition 6.1.}

  We say that the group $G$ is Noetherian
if and only if for every set of equations $E$, there exists a finite subset $E_0$ of $E$ such that
$V_G(I_E)=V_G(I_{E_0})$.

\bigskip

 The problem of the existence of a solution of a finite set of equations defined on a group $G$
 has been studied by many authors.
This situation is analog the study of the roots of polynomial equations in field theory. In fact, Scott [7] has defined
a notion of algebraically closed group that we recall.

\medskip

{\bf Definition 6.2.}

A set $E$ of equations and inequations defined on a group $G$ is consistent, if there exists an embedding $G\rightarrow H$ such that $E$ has a solution
in $H$.

\medskip

An example of a set of equations and inequations which is not consistent is: $x^2 =1, x^3 =1,x\neq 1$.

\medskip

{\bf Definitions 6.3.}

A group is algebraically closed if and only if every finite set of equations consistent defined on $G$ has a solution in $G$.

Scott has shown that every group can be imbedded in an algebraically closed group.
We show now the following proposition (See also Baumslag and al [2] p. 62)

\medskip

{\bf Proposition 6.1.}

{\it Let $G$ be an algebraically closed group, then for every finitely generated maximal $M$ ideal of $G[X_1,...,X_n]$,
there exists elements $a_1,...,a_n$ in $G$ such that $M$ is generated by $\{X_1a_1^{-1},...,X_na_n^{-1}\}$.}

\medskip

{\bf Proof.}

Let $H$ be the quotient $H=G[X_1,...,X_n]/M$, and $x_1,...,x_n$, the image of $X_1,...,X_n$ by the natural
projection $G[X_1,...,X_n]\rightarrow H$. Then $(x_1,...x_n)$ is a solution of every equation defined by
elements of $M$. Since $G$ is algebraically closed and $M$ is finitely generated, there exists an element $(a_1,...,a_n)$ in $G^n$
which is a solution of the equations defined by elements of $M$. This implies that the ideal generated by $X_1a_1^{-1},...,X_na_n^{-1}$
 contains $M$ henceforth this ideal   is equal to $M$ since $M$ is maximal.

\bigskip

{\bf Definition 6.4.}

A topological set is Noetherian if and only if every descending chain of closed subspaces
$Y_1\supseteq Y_2....\supseteq Y_i....$ stabilizes. This means that there exists $n$ such that $Y_i=Y_n$
for every $n\geq i$.

Let $H$ be an object of $C(G)$,  $H$ is a Noetherian group if and only if every
 family of ideals $I_1\subset I_2\subset...\subset I_n...$ stabilizes, this is equivalent to saying
 that there exists an $i$ such that for every $n\geq i$, $I_n = I_i$.

 Remark that $Spec_G(H)$ is Noetherian if and only if $H$ is a Noetherian group.

 \medskip

 {\bf Proposition 6.2.}

 {\it The object $H$ of $C(G)$ is Noetherian if and only if every ideal of $H$ is finitely generated.}

 \medskip

 {\bf Proof.}

 Suppose that $H$ is Noetherian, let $I$ be an ideal of $H$, consider the set $A$ whose objects are
finitely generated  ideals contained in $I$ ordered by the inclusion. Every totally ordered family $(I_l)_{l\in L}$ of $A$
 has a maximal element. If not we can construct a subfamily $(I_n)_{n\in N}$ of $(I_l)_{l\in L}$
 such that $I_n$ is strictly contained in $I_{n+1}$. This is in contradiction with the fact that
 $H$ is Noetherian. Let $M$ be a maximal element of $A$, $M= I$, if not there exists an element $x$
 of $I$ which is not in $M$ and the ideal generated by $M$ and $x$ contains strictly $M$.

 Suppose that every ideal of $H$ is finitely generated. Let $I_0\subset I_2...\subset I_n...$ be
 an ascending sequence of ideal, $\bigcup_{n\in N}I_n = I$ is a finitely generated ideal by
 $x_1,...,x_p$. The element $x_i$ belongs to $I_{m(i)}$. Let $m$ be the maximum of $m(i)$, $I = I_m$.

\medskip

{\bf Proposition 6.3.}

{\it Suppose that $Spec_G(H)$ is a Noetherian topological space, then the radical of $H$ is the intersection
of a finite number of primes ideals. In particular if the intersection of all the elements of $Spec_G(H)$ is the
neutral element, then there exists a finite number of elements of $Spec_G(H)$ whose intersection is the neutral element.}

\medskip

{\bf Proof.}

We know that a Noetherian topological space is the union of a finite number of irreducible closed subsets.
This implies that there exists a finite number of primes ideals $P_1,...,P_n$ such that $Spec_G(H)=\bigcup_{i=1,...,n}V(P_i)$.
The proposition 4.1 implies that $Spec_G(H) = V(\bigcap_{i=1,...,n}P_i)$. This implies that every prime ideal $P$ of $H$
contains $\bigcap_{i=1,...,n}P_i$, thus $\bigcap_{i=1,...,n}\subset Rad_G(H)$. It results that $Rad_G(H) =\bigcap_{i=1,...,n}P_i$
since $Rad_G(H)$ is the intersection of all the elements of $Spec_G(H)$.

\bigskip

{\bf 7. Dimension of affine schemes.}

\medskip

We want to introduce the notion of dimension of an affine schemes. In commutative algebra, the dimension
of a ring can be defined to be its height that is, the maximal length of  chains $0\subset P_1...\subset P_n$ of prime ideals
such that $P_i$ is distinct of $P_{i+1}$. We can adopt the same definition here for the prime ideals of
elements of $H$ for every object $H$ of $C(G)$. But the following proposition shows that the height of $G*Z$
can be infinite.

\medskip

{\bf Proposition 7.1.}

{\it Let $G$ be a group which does not have divisors of zero, then the height of $G*Z$ is infinite.}

\medskip

{\bf Proof.}

Let $I_n$ be the ideal of $G*Z$ generated by $2^n$, $(G*Z)/I_n = G*(Z/2^n)$. Suppose that $G*Z/2^n$
have divisors of zero $x,y$. For every $g\in G$, we have $[x,y^g] = 1$, $x$ and $y$ cannot be in a subgroup conjugated to
$G*1$ since $G$ does not have zero divisors; $x,y$ cannot be in a subgroup conjugated to $1*Z$ in this case, $[x,y^g]$ is
not trivial for every non trivial element of $g\in G$.
 If $x$ is not in a conjugated of $G*1$ and $1*Z$, then [5] page 187 shows that for every $g\in G$, $y^g$ and $x$ are the
power of the same element; this is impossible. We conclude that $I_n$ is a prime, thus the sequence $I_n\subset I_{n-1}...\subset I_0$
is an infinite sequence of distinct prime ideals.

\medskip

In classical algebraic geometry, the dimension of a scheme is also the the transcendence dimension
of its field of fraction. The corresponding notion here is the $G$-rank:

\medskip

{\bf Definition 7.1.}

Let $H$ be an object of $C(G)$, a family of $G$-generators of $H$ is a family of elements $(x_i)_{i\in I}$ of $H$
such that $G$ and $(x_i)_{i\in I}$ generate $H$. The $G$-rank of $H$ is the minimal cardinal of the set $I$ such that there
exists a family of generators $(x_i)_{i\in I}$ of $H$.

\medskip

 We have the following result:

\medskip

{\bf Proposition 7.2.}

{\it Let $G$ be a simple group, and $G^n$ the product of $n$ copies of $G$
endowed with the $G$-structure defined by the embedding of $G$ as the first factor, $(G,1,...,1)$ then $Spec_G(H)$ contains
only one element the prime $(1,G,...,G)$ and the $rank_G(G^n) = n-1$.}

\medskip

{\bf Proof.}

Let $I$ be an ideal of $G^n$, suppose that $I$ contains an element $u =(g_1,...,g_n)$ such that $g_1$ is different of
the neutral element of $G$. For every $g\in G$, let $i_g = (g,1,...,1)$, we have:
 $i_gui_g^{-1}= (gg_1g^{-1},g_2,...,g_n)\in I$. Since $G$ is simple, there exists $g\in G$ such that $gg_1g^{-1}\neq g_1$.
 This implies that $u^{-1}i_gui_g^{-1}$ is a non trivial element of $G$. This is a contradiction since $I$ is an ideal, thus
 $I$ is contained in $\{1\}\times G^{n-1}$. If $I$ is distinct of $\{1\}\times G^{n-1}$, there exists an element $v\in \{1\}\times G^{n-1}$ who does
 not belong to $I$. Let $p:G^n\rightarrow G^n/I$ be the natural projection, for every element $g\in G$, we have
 $[g,v] =1$ thus $[p(g),p(v)] =1$. This implies that $G^n/I$ is not a domain, thus $I$ is not a prime.

 \medskip

{\bf Remark.}

Suppose that the center of $G$ is not trivial, then for every element of $C(G)$, $Spec_G(H)$ is empty and the $G$-rank
of $H$ needs not to be $0$ or $-\infty$.

\medskip

{\bf Definition 7.2.}

Let $H$ be an object of $C(G)$:

-  if $Spec_G(H)$ is empty, the dimension of $Spec_G(H)$ is $-\infty$.

if $Spec_G(H)$ is not empty,  the dimension of $Spec_G(H)$ is the $G$-rank of $H/Rad_G(H)$
the quotient of $H$ by the intersection $Rad_G(H)$ of all the elements of $Spec_G(H)$.

\medskip

{\bf Remarks.}

If $H$ is a $G$-domain, then $Rad_G(H) = \{1\}$ and the dimension of $Spec_G(H)$ is $rank_G(H)$.

Let $G$ be a group without zero divisors, and $H = G[x_1,...,x_n]$, $H$ is a $G$-domain and the dimension of $H$ is $n$.

If $H$ is a $G$-domain, then dimension of $Spec_G(H) = 0$ is equivalent to saying that $H$ is isomorphic
to $G$.

Let $H = G^n$ endowed with the $G$-structure defined by the imbedding $(G,1,...,1)$ of $G$ in $H$, then
$Rad_G(H) = (1,G,...,G)$ and $H/Rad_G(H) \simeq G$; this implies that the dimension of $Spec_G(H) =0$.

\medskip

We know that a if $H$ is $G$-simple, then $Spec_G(H)$ contains only one element, but the $G$-rank of $H$ is not
necessarily $1$. But if $G$ is algebraically closed, we have the following result:

\medskip

{\bf Proposition.}

{\it Suppose that $G$ is algebraically closed and $H$ a $G$-domain of finite type which is $G$-simple, then $H=G$, thus $dim_G(H) = 0$.}

\medskip

{\bf Proof.}

There exists a positive integer $n$ and an ideal $I$ of $G[x_1,...,x_n]$ such that $H = G[x_1,...,x_n]/I$.
Since $G$ is algebraically closed, there exists a maximal ideal $M$ (defining by a solution of the equations defined by a
 finite set of generators of $I$ containing $I$) such that $G[x_1,...,x_n]/M = G$. Let $M'$ be the image of $M$
 by the projection map $G[x_1,...,x_n]\rightarrow H$, $M'$ is an ideal of $H$, since $H$ is $G$-simple, we conclude
 that $M'=\{1\}$ thus $M=I$. We deduce that $H = G$ since $ H = G[x_1,...,x_n]/M$.

\medskip

Let $H$ be a group, the co-rank of $n$ is the maximal integer $n$ such that there exists a surjective
morphism: $H\rightarrow F[x_1,...,x_n]$. We say that the $G$-co-rank of the object $n$ in $C(G)$ is $n$ if $n$ is the
maximal integer such that there exists a $G$-morphism: $H\rightarrow G[x_1,...,x_n]$.

Contrary to the rank, a class of ideals allows to express geometrically the idea of the co-rank,

\medskip

{\bf Definition.}

Let $H$ be an element of $C(G)$, a $G$-free ideal $I$ of $H$ is an ideal such that $H/I =G[x_1,...,x_n]$.

The co-dimension of $Spec_G(H)$ is the maximum integer $n$ such that there exists an ideal $I$ such that $H/I = G[x_1,...,x_n]$.

\medskip

{\bf Proposition 7.4.}

{\it Let  $H$ be an element of $C(G)$, the co-dimension of $Spec_G(H)$ is the length of the maximal
sequence $P_0\subset P_2\subset...\subset P_n$ such that for $i=1,...,n$ $P_i$ is $G$-free and  $P_{i+1}$ contains strictly $P_i$.}

\medskip

{\bf Proof.}

Suppose that the co-dimension of $Spec_G(H)$ is $n$, there exists a surjective $G$-morphism
$f:H\rightarrow G[x_1,...,x_n]$. We denote by $P_0$ the ideal$f^{-1}(1)$ of $H$. Let $I_l$ be the normal subgroup
 of $G[x_1,...,x_n]$ generated by $x_1,...,x_l$, and for $i=1,...,n$, $P_i =f^{-1}(I_i)$.
The family of ideal $(P_i)_{i\in\{0,...,n\}}$ is a chain of free ideals of length $n$ such that $P_i$ is strictly included in $P_{i+1}$.

Let $P_0\subset P_2\subset...\subset P_n$ be a maximal chain of $G$-free ideals such that  $P_{i+1}$ contains strictly
 $P_i$. Remark that the $G$-co-rank of $H$ is inferior or equal to $n$, since
if the co-rank is superior to $n$ the first stage of the proof shows that we can construct a chain of $G$-free prime ideals of
length strictly greater than $n$ and whose $G$-co-rank is strictly decreasing. Thus it is enough to show
that $G/P_0\simeq G[x_1,...,x_n]$.

Suppose that the $G$-rank of $G/P_0$ is strictly inferior to $n$; this implies that
there exists, $i\leq n$ such that co-rank$(P_i)$ =co-rank$(P_{i+1})$. Let $l_i: H\rightarrow H/P_i$ be the
natural projection, $H/P_{i+1}$ is the quotient of $H/P_i$ by $l_i(P_{i+1})$. Write $H/P_i =G[x_1,...x_{n_i}]$ and
$H/P_{i+1} = G[y_1,...,y_{n_i}]$. The quotient of $G[x_1,...,x_{n_i}]$ by the normal subgroup generated by $G$ is
isomorphic to the free group $F[x_1,...,x_{n_i}]$. The quotient of $G[y_1,...,y_{n_i}]$ by the normal subgroup generated
by $G$ is isomorphic to the quotient of $F[x_1,...,x_{n_i}]$ by the normal subgroup generated by  the image of $l_i(P_{i+1})$
in $F[x_1,...,x_{n_i}]$. Since this group is isomorphic to the free group of rank $n_i$ $F[y_1,...,y_{n_i}]$, we  deduce
that $l_i(P_{i+1})$ is included in the normal subgroup of $G[x_1,...,x_{n_i}]$ generated by $G$. This in
contradiction with the fact that the quotient of $G[x_1,...,x_{n_{i+1}}]$ by the image of $l_i(P_{i+1})$ is a free product of $G$
and a free group.

\medskip

The $G$-free ideals of $H$ generate a topology as shows the following proposition:

\medskip

{\bf Proposition 7.5.}

{\it Let $G$ be a group which does no have divisors of zero, and  Let $H$ be an element of                                                                  $C(G)$, for every normal subgroup $I$ of $H$, denote $V_F(I) = \{P:$ $P$ is a $G$-free prime ideal
and $I\subset P\}$, then for every normal subgroups $I,J$ of $H$, $V_F(I)\bigcup V_F(J) = V_F([I,J])$. Let
$(I_a)_{a\in A}$ be a family of normal subgroups of $H$. Denote by $I_A$ the normal subgroup  generated by $(I_a)_{a\in A}$,
then $V_F(I_A) =\bigcap_{a\in A}V_F(I_a)$.}

\medskip

{\bf Proof.}

Firstly we show that $V_F(I)\bigcup V_F(J) = V_F([I,J])$. Let $P$ be an element of $V_F(I)\bigcup V_F(J)$, $P$
contains $I$ or $J$ this implies that $P$ contains $[I,J]$. Let $P$ be an element of $V_F([I,J])$. Suppose
that $P$ is neither an element of $V_F(I)$ nor an element of $V_F(J)$. This implies that there exists an element $x\in I$ and $y\in J$ which
are not in $P$. For every $g\in G$, $[x,y^g]\in P$. Let $p:H\rightarrow H/P$ be the natural projection,
we have $[p(x),p(y)^g] = 1$. Since $H/P$ is $G$-free, we deduce that $x$ and $y$ are elements of $G$, this is in contradiction
with the fact that $G$ does not have divisors of zeros.

Now we show that $V(I_A) = \bigcap_{a\in A}V(I_a)$. Let $P$ be an element of $V(I_A)$, since $I_a\subset I_A \forall a\in A,$
$P\subset V(I_a) \forall a\in A$. Let $P$ be an element of $\bigcap_{a\in A}V(I_a)$, $\forall a\in A,I_a\subset P$,
this implies that $I_A\subset P$ and $P\in V(I_A)$.

\medskip

{\bf Remark.}

Let $G$ be a group which does not have zero divisors, and $H$ an object of $C(G)$, the previous proposition
shows that the family of subspaces $V_F(I)$ where $I$ is any normal subgroup in $H$ defines a topology on
the space $Spec_G^F(H)$ of $G$-free prime ideals of $H$. But if $f:H\rightarrow H'$ is a morphism between
the elements $H$ and $H'$ of $G$ and $P$ is an element of $Spec_G^F(H)$, $f^{-1}(P)$ is not necessarily a $G$-free
prime of $H$. For example let $G = PSL(n,Z)$ and  $H' = PSL(n,Z)[X]$. Let $H$ be the subgroup of $H'$ generated by $PSL(n,Z)$ and $XNX^{-1}$
where $N$ is a nilpotent subgroup of $PSL(n,Z)$. The canonical imbedding $i:H\rightarrow H'$ is a $G$-morphism.
The normal trivial subgroup $I=\{1\}$ of $H'$  is a prime ideal of $H$, but the quotient of $H$ by $i^{-1}(I)$
is not $G$-free, since $H$ is not $G$-free.

The theorem of Kurosh  says that  a subgroup of a free product of groups $H_1*H_2$ is the free product of a free group
and groups isomorphic to subgroups of $H_1$ and $H_2$. Thus if we choose $G=F_n n\geq 2$, a $G$-morphism $f:H\rightarrow H'$
induces a morphism $Spec_G^F(H')\rightarrow Spec_G^F(H)$ of $G$-free spectra. In fact in this case the topology
is defined by a Top couple (see definition 3.1).

\medskip

{\bf Definition 7.4.}

Let $G$ be a group without divisors of zero, and $H$ an element of $C(G)$; an element
 $P$ of $Spec_G(H)$ is a $G$-point if and only if $H/P = G$.

\medskip

{\bf Remark.}

If $P$ is a $G$-point, the restriction of the projection $p:H\rightarrow H/P$ to $G$ is an isomorphism, thus
this extension splits.

\medskip

{\bf Proposition 7.6.}

{\it let $H$ be a be a $G$-domain if $G$ is a normal subgroup of $H$, for every prime $P$ of $H$, we have
of $[G,P] = 1$. If there exists a $G$-point $P$ in $H$, then $H = P\times G$.}

\medskip

{\bf Proof.}

Let $P$ be a prime ideal.  We have $[P,G]\subset P\cap G =\{1\}$ since $P$ is an ideal and $G$ is a normal subgroup of $H$.

If $P$ is a $G$-point of $H$, then $H$ is an extension of $G$ by $P$ which splits, since $[P,G] = 1$ we deduce the result.

\bigskip

{\bf 8. Galois theory.}

\medskip

Let $G$ be a group, a polynomial of $G$ is an element of $G[X]$.
The Galois theory in classical field theory enables to study the solutions of polynomial functions defined
on a field. One of the property essential in that study is the existence of factorization: let $P$ be a
polynomial defined on a field, and $a$ a root of $P$, we can write $P = (X-a)Q$; if $deg(P)$ is the degree of $P$, we have:
 $deg(Q) = deg(P) -1$. Here we do not have factorization system, in fact we can have a polynomial
 which has an infinite numbers of roots: Let $G = Gl(2,C)$ be the group of automorphisms of the
 complex plan;  the solutions of the equation defined on $G$ by $X^2 = 1$ is the set of symmetries of
 the complex plan. We can introduce anyway some notions related to Galois theory.

 \medskip

 {\bf Definitions 8.1.}

 Let $G$ be a group, and $H$ an element  of $C(G)$, an element $a\in H$ is algebraic over $G$ if $a$ is
 the root of a polynomial whose coefficients are in $G$.

 An element $H$ of $C(G)$ is an algebraic extension of $G$ if and only if every element of $H$ is algebraic
 over $G$.

  We say that an algebraic extension of $G$ is finite if and only if the action of $G$ on $H$ by left multiplications
  has a finite number of orbits

 \medskip

 {\bf Proposition 8.1.}

{\it Suppose that $G = \{1\}$, then a group $H$ is an algebraic extension of $G$ if and only if every
element of $H$ has a finite order. In particular a group is a finite extension of $\{1\}$ if and only if
it is a finite group.}

\medskip

{\bf Proof.}

If $G=\{1\}$, then $G[X] = F[X]= Z$ and the polynomials of $\{1\}$ are $X^n$. A group $H$ is algebraic over $\{1\}$
if and only if for every $h$ in $H$, there exists an element $P$ of $F[X]$ such that $P(h) =1$. This is equivalent
to the existence of an integer $n$ such that $h^n =1$.

Suppose that $H$ is a finite group. Let $n_H$ be the order of $H$, for every element $h\in H$, we have $h^{n_H}=1$;
this implies that $H$ is algebraic over $\{1\}$. The left action of $\{1\}$  has finite orbit since $H$ is finite.
Conversely, if the left action of $\{1\}$ on $H$ has finite orbits, then $H$ is a finite group.

\medskip

{\bf Remark.}

An algebraic extension  $H$ of $\{1\}$ is not necessarily finite even if $H$ is finitely generated. In 1904 Burnside asks the question wether every
finitely generated group for which every element has a finite order is finite. There are groups which have been constructed
answering negatively the problem of Burnside.

\medskip

{\bf Proposition 8.2.}

{\it Let $H$ be an element of $C(G)$, suppose that $H$ is the union of a finite number of orbits for
 the action of $G$ on $H$ defined by the left multiplication, then $H$ is algebraic over $G$. In particular if $G$ is
a normal subgroup of $H$ such that $G/H$ is finite, $H$ is algebraic over $G$.}

\medskip

{\bf Proof.}

Let $x$ be an element of $H$; if  $H$ is the union of a finite number of orbits for the left action of $G$, there exists
distinct integers $n$ and $m$ with $n<m$  such that $x^n$ and $x^m$ are in the same orbit. This is equivalent
to saying that there exists $h\in H$, $g_1,g_2\in G$ such that $x^n = g_1h$ and $x^m =g_2h$, we deduce that
$h={g_1}^{-1}x^n={g_2}^{-1}x^m$. It results that $g_1{g_2}^{-1}x^{m-n}=1$. Thus $x$ is algebraic over $G$.

\medskip

{\bf Proposition 8.3.}

{\it Let $H$ be an algebraic extension of $G$, then the $G$-co-rank of $H$ is strictly inferior to $1$.}

\medskip

{\bf Proof.}

Suppose that the $G$-co-rank of $H$ is superior or equal to $1$. There exists a $G$-surjection:
$s:H\rightarrow G[x_1,...,x_n]$. Let $y_1,...,y_n$ a family of algebraic elements of $H$ which generate $H$.
There exists $i\in \{1,...,n\}$ such that $s(y_i)$ is not in $G$. Since $y_i$ is algebraic, there exists a $G$-polynomial
$P$ such that $P(y_i) =1$. We also have $P(s(y_i)) = 1$. This is impossible since  $s(y_i)$ is a  non trivial words
written with elements of $G$ and with elements of the free group generated by $x_1,...,x_n$.

 \medskip

{\bf Proposition 8.4.}

{\it Let $H$ be an element of $C(G)$ and $K$ an element of $C(H)$, suppose that $H$ is a finite algebraic
extension of $G$ and $K$ is a finite algebraic extension of $H$; then $K$ is a finite algebraic extension of $G$.}

\medskip

{\bf Proof.}

To show this proposition it is enough to show that the action of $G$ on $K$ has a finite number of orbits
by applying proposition 8.2.

Since $H$ is a finite algebraic extension of $G$ and $K$ a finite algebraic extension of $H$, we know that  the orbit spaces $H/G$
and $K/H$ are finite; we denote
by $h_1,...,h_m$  elements of $H$ such that the union of $Gh_1,...,Gh_m$ is $H$,
 and by $k_1,...,k_n$  elements of $K$ such that the union of $Hk_1,...,Hk_n$ is $K$. Let $x$ be an element of $K$,
there exists $1\leq i\leq n$ such that $x\in Hk_i$, since $H$ is the union of $Gh_1,...,Gh_m$, we deduce that
$x\in \{Gh_1k_i,...,Gh_mk_i\}$. This is equivalent to saying that the union of $Gh_jk_i$  $1\leq j\leq m$, $1\leq i\leq n$ is $K$.

\medskip

 {\bf Definition 8.2.}

 Let $H$ be an algebraic extension of $G$, the Galois group $Gal_G(H)$ of $H$ is the group of automorphisms
 of $H$ whose restriction to $G$ is the identity.

\bigskip
\bigskip

{\bf References.}

\bigskip

1.  Amaglobeli. M.G Algebraic sets and coordinate groups for a free nilpotent group of nilpotency class 2. Sibirsk. Mat. Zh. Volume 48 p. 5-10.

\smallskip

2.  Baumslag, G,  Miasnikov, A.  Remeslennikov, V.N. Algebraic geometry over groups I. Algebraic sets and ideal theory. J. Algebra. 1999, 219, 16–79.

\smallskip

3. M. Chabashvili. Lattice isomorphisms of nilpotent of class 2 Hall W powers groups. Bulletin of the Georgian National Academic of Sciences, volume3
2 2008

\smallskip

4. A.Grothendieck, \'El'ements de g\'eom\'etrie alg\'ebrique I.Publications math\'ematiques de l'I.H.E.S 4, 5-228

\smallskip

5. A. Karass, W. Magnus, D. Solitar, combinatorial group theory, Dover publication 1976

\smallskip

6. Micali, A. Alg\`ebres g\'en\'etiques, Cahiers de math\'ematiques de Montpellier 1985

\smallskip

7. Scott W.R Algebraically closed groups Proc. Amer. Math. Soc. 2 (1951) 118-121

\smallskip

8. Serre J-P. G\'eom\'etrie alg\'ebrique, g\'eom\'etrie analytique, Annales de l'Institut Fourier, Grenoble
 t. 6 1955-1956 1-42.

 \smallskip

9. J-P. Serre, Arbres,  amalgames, $SL_2$ Ast\'erisque 46 1977.

\smallskip

10. Jacques Tits; Sous-alg\`ebres des alg\`ebres de Lie semi-simples. S\'eminaire Bourbaki 1954-1956 p. 197-214

\smallskip

11. Tits, J; Free groups in linear groups, Journal of Algebra 20 (2) 250-270 (19720

\smallskip

12. Tsemo, A. Scheme theory for groups and Lie algebra, International Journal of Algebra 5. 2011 139-148

\smallskip

13. Tsemo, A. Algebraic geometry on groups: the theory of curves In preparation.

\end{document}